\documentclass[11pt,reqno]{amsart}
\usepackage{euscript,amssymb}
\usepackage[arrow,matrix,curve]{xy}
\usepackage{array,longtable}
\usepackage{a4wide}

\newenvironment{m-theorem}{%
\vskip5pt\refstepcounter{stff}\trivlist \itemindent 0pt
\item[\hskip\labelsep\bf Theorem \thestff]%
\it\ignorespaces}{\endtrivlist\vskip5pt}%
\newenvironment{m-proposition}{%
\vskip5pt\refstepcounter{stff}\trivlist \itemindent 0pt
\item[\hskip\labelsep\bf Proposition \thestff]%
\it\ignorespaces}{\endtrivlist\vskip5pt}%
\newenvironment{m-corollary}{%
\vskip5pt\refstepcounter{stff}\trivlist \itemindent 0pt
\item[\hskip\labelsep\bf Corollary \thestff]%
\it\ignorespaces}{\endtrivlist\vskip5pt}%
\newenvironment{m-lemma}{%
\vskip5pt\refstepcounter{stff}\trivlist \itemindent 0pt
\item[\hskip\labelsep\bf Lemma \thestff]%
\it\ignorespaces}{\endtrivlist\vskip5pt}%
\newenvironment{m-definition}{%
\vskip5pt\refstepcounter{stff}\trivlist \itemindent 0pt
\item[\hskip\labelsep\bf Definition \thestff]%
\ignorespaces}{\endtrivlist\vskip5pt}%
\newenvironment{m-notation}{%
\vskip5pt\refstepcounter{stff}\trivlist \itemindent 0pt
\item[\hskip\labelsep\bf Notation \thestff]%
\ignorespaces}{\endtrivlist\vskip5pt}%
\newenvironment{m-example}{%
\vskip5pt\refstepcounter{stff}\trivlist \itemindent 0pt
\item[\hskip\labelsep\bf Example \thestff]%
\ignorespaces}{\endtrivlist\vskip5pt}
\newenvironment{m-remark}{%
\vskip5pt\refstepcounter{stff}\trivlist \itemindent 0pt
\item[\hskip\labelsep\bf Remark \thestff]%
\ignorespaces}{\endtrivlist\vskip5pt}
\newenvironment{m-question}{%
\vskip5pt\refstepcounter{stff}\trivlist \itemindent 0pt
\item[\hskip\labelsep\bf Question.]%
\ignorespaces}{\endtrivlist\vskip5pt}%
\newenvironment{thm-nono}{%theorem no-number
\vskip5pt\trivlist \itemindent 0pt
\item[\hskip\labelsep\bf Theorem.]%
\it\ignorespaces}{\endtrivlist\vskip5pt}%
\newenvironment{m-thank}{%
\vskip5pt\trivlist \itemindent 0pt
\item[\hskip\labelsep\it Acknowledgments]%
\ignorespaces}{\endtrivlist\vskip5pt}%

\let\mt\mapsto

\font\tenmsa=msam10 %
\newcommand\hdashpiece{%
{\vrule height2.75pt depth-2.35pt width2.3pt \kern1.7pt}}%
\newcommand\hdashpieces{%
{\hdashpiece\hdashpiece\hdashpiece\hdashpiece}}%

\newcommand\dashar{\mathrel{%
\hdashpieces\kern-0.4pt\hbox{\tenmsa K}}}%

\let\euf\EuScript %use of package ``euscript'' required

\let\mbb\mathbb

\DeclareFontFamily{OT1}{rsfs}{}
\DeclareFontShape{OT1}{rsfs}{n}{it}{<->rsfs10}{}
\DeclareMathAlphabet{\crl}{OT1}{rsfs}{n}{it}

\let\ovl\overline

\let\nit\noindent

\let\vphi\varphi

\let\lan\langle
\let\ran\rangle
\newcommand\lran[1]{{\lan #1\ran}}

\newcommand\Aut{\operatorname{\textrm{Aut}\kern1pt}}
\newcommand\cAut{\operatorname{\mathcal{A}\kern-1pt\textit{ut}\kern1pt}}
\newcommand\End{\operatorname{\rm{End}\kern1pt}}
\newcommand\cEnd{\operatorname{\mathcal{E}\kern-1pt\textit{nd}\kern1pt}}

\newcommand\cHom{\operatorname{\mathcal{H}\kern-1pt\textit{om}\kern1pt}}

\newcommand\Pic{\mathop{\rm Pic}\nolimits}

\newcommand\invq{{\slash\kern-2.5pt\slash}}

\newcommand\rk{{\rm rk}}

\newcommand\Grs{{\rm Gr}}

\numberwithin{equation}{section}
\numberwithin{figure}{section} 

\let\L\Lambda

\newcommand\bbC{{\mbb C}}

\newcommand\bbk{\mbox{\rm I\kern-1.5pt k}}
\newcommand\sbbk{\hbox{\scriptsize I{\kern-.8pt}k}}

\newcommand\bone{{1\kern-0.57ex\rm l}}

\newcommand\eL{{\euf L}}

\newcommand\eO{{\euf O}}

\newcommand\codim{{\rm codim}}

\newcommand\SO{{\rm SO}}
\newcommand\Sl{{\rm SL}}
\newcommand\Sp{{\rm Sp}}

\let\ges\geqslant

\newcommand\res{\mathop{\rm res}\nolimits}

\newcommand\DD{{D}}
\newcommand\LL{{\ell}}

\newcommand{\lth}{{\rm length}}
\newcommand{\LG}{\mathop{LG}\nolimits}

\author{Mihai Halic}
\address{}

\begin{document}

\title{Splitting criteria for vector bundles\\ on minuscule homogeneous varieties}

\begin{abstract}
I prove that a vector bundle on a minuscule homogeneous variety splits into a direct sum of line bundles if and only if its restriction to the union of two-dimensional Schubert subvarieties splits. A case-by-case analysis is done. 
\end{abstract}

\maketitle\markboth{\sc Mihai Halic}{\sc Splitting of vector bundles on minuscule homogeneous varieties}

%%%%%%%%%%%%%%%%%%%%%%%%%%%%%%%%%%%%%%%%%%%%
%%%%%%%%%%%%%%%%%%%%%%%%%%%%%%%%%%%%%%%%%%%%

\section*{Introduction}

As stated in the abstract, the goal of this note is to prove the following: 
\begin{thm-nono} 
A vector bundle on a \emph{minuscule} homogeneous variety $X$ with cyclic Picard group splits if and only if its restriction to the union of $2$-dimensional Schubert subvarieties of $X$ splits. 
\end{thm-nono}
The list of minuscule flag varieties can be found e.g. \cite[Ch. 5, \S2]{hiller-book}. For the groups of type $A_n,B_n,D_n$ one gets respectively the Grassmannians, the even dimensional quadrics, and the spinor varieties. There are two `exceptional' cases, the Cayley plane and the Freudenthal variety, corresponding to the groups $E_6,E_7$ respectively. 

A case-by-case analysis shows that the statement above is very down-to-earth.   Minuscule homogeneous varieties have only one irreducible $2$-dimensional Schubert subvarieties, isomorphic to $\mbb{CP}^2$, with two exceptions: the Grassmannians $\Grs(k,n)$, $1<k<n-1$, and the $4$-dimensional quadric, which have two irreducible $2$-dimensional Schubert varieties, both isomorphic to $\mbb{CP}^2$, intersecting along a Schubert line.

The problem of deciding the splitting of vector bundles goes a long way back. One of the first and widely known criterion is that of Horrocks \cite{horr}, saying that a vector bundle on the projective space splits if and only if it does so along a $2$-plane. Cohomological splitting criteria for vector bundles on Grassmannians and quadrics have been obtained in \cite{ottv,mals}. Recent progress is achieved in \cite{mos}, where the authors proved a splitting criterion for uniform vector bundles of \emph{low rank} on more general Fano varieties with cyclic Picard group. (The term `uniform' refers to the fact that the splitting type of the vector bundle is the same along \emph{all} the rational curves representing a fixed homology class. Earlier articles which investigated uniform vector bundles on projective spaces and Grassmannians are \cite{vdv,ehs,bal,guy}.) One of the main applications yields splitting criteria for uniform vector bundles of low rank on Hermitian symmetric spaces, also known as \emph{co-minuscule} homogeneous varieties. The assumption on the rank of the vector bundles to be sufficiently small is essential: there are easy examples of homogeneous, thus uniform, non-split vector bundles; optimal upper bounds for the rank are given in \textit{op. cit.}

The theorem above requires to probe the splitting of vector bundles along a single $2$-plane or a wedge of two $2$-planes; moreover, the \emph{rank} is \emph{arbitrary}. Furthermore, there is an elementary cohomological test (cf. proposition \ref{prop:x2}) to verify whether a vector bundle on the projective plane splits or not. From a `computational' point of view, these facts are an advantage because, unless additional information is available, it is difficult to check whether a vector bundle is uniform or not. Indeed, the generic splitting type of any vector bundle being constant, one has to verify the uniformity on every line 
(within a suitable homology class). 

One of the essential tools used for proving the main result is Pieri's formula in Schubert calculus, which is recalled in the first section. The main result itself is proved in section \ref{sct:result}, and detailed in section \ref{sct:expl}.

%%%%%%%%%%%%%%%%%%%%%%%%%%%%%%%%%%%%%%%%%%%%
%%%%%%%%%%%%%%%%%%%%%%%%%%%%%%%%%%%%%%%%%%%%

\section{Preliminaries}\label{sct:frame}

Here we recall from \cite[\S1]{lms}, \cite[Ch.~IV,V]{hiller-book} a few basic facts about the Bruhat decomposition of projective homogeneous varieties and their Schubert subvarieties; a thorough description of the cohomology ring can be found in \cite{bgg,dem}. Let $G$ be a connected, reductive linear algebraic group and $P\subset G$ a parabolic subgroup; denote by $W_G$ and $W_P$ the Weyl groups of $G$ and of the Levi subgroup of $P$ respectively. The Bruhat decomposition states that 
\begin{equation}\label{eq:bruhat}
X:=G/P=\underset{w\in W_G/W_P}{\mbox{$\coprod$}}BwP. 
\end{equation}
Each $W_P$-orbit in $W_G$ admits a unique representative of shortest length, and the dimension of the corresponding cell $BwP$ equals $\lth(w)$; let $W^P$ be the set of shortest representatives. We denote by $w_0\in W_G$, $w_{0,P}\in W_P$, $w_0^P\in W^P$ be the longest words; they are related by $w_0=w_0^Pw_{0,P}$.  By definition, the Schubert subvarieties of $G/P$ are the closures of the Bruhat cells: 
$$
\ovl{BwP},\quad\dim \ovl{BwP}=\lth(w),\quad\forall\,w\in W^P.
$$ 
Since we are interested in intersection products, it is more convenient to index the Schubert subvarieties by their codimension. The Poincar\'e duality on $G/P$ acts as: 
\begin{equation}\label{eq:dual}
W^P\to W^P,\quad w\mt \check w:=w_0ww_{0,P}.
\end{equation}
(The right-hand-side is not the reduced decomposition; it can be simplified by using the commutation relations in $W_G$.) We denote 
\begin{equation}\label{eq:Xw}
X(w):=\ovl{B\check{w}P},\quad\codim X(w)=\lth(w),\quad\forall\,w\in W^P.
\end{equation}
The (reversed) \emph{Bruhat order} on $W^P$ is defined by: 
$$
w\prec w'\;\Leftrightarrow\;X({w'})\subsetneq X({w}).  
$$
It defines an oriented graph (the Hasse diagram) with vertices $W^P$, and edges 
$$
\{(w,w')\mid w\prec w',\;\lth(w')=\lth(w)+1\}.\quad(\text{Notation:}\;w\to w').
$$ 
The Bruhat cells are isomorphic to affine spaces, so \eqref{eq:bruhat} is a cell-decomposition of $X$, and therefore $\{X(w)\}_{w\in W^P}$ is a $\mbb Z$-basis of the cohomology ring $H^*(X;\mbb Z)$ of $X$; in particular, $H^{\rm odd}(X;\mbb Z)=0$. The Picard group is isomorphic to $H^2(X;\mbb Z)$, generated by the classes of the divisors $\DD_\alpha:=X(\check{\tau}_\alpha)$, where $\alpha$ is a simple root such that $\tau_\alpha\in W^P$. (Here $\tau_\alpha$ stands for the reflection defined by $\alpha$.) The intersection product between $\DD_\alpha$  and the Schubert variety $X(w)$ is given by Pieri's formula: 
\begin{equation}\label{eq:pieri}
\big[\DD_\alpha\cdot X(w)\big]=\biggl[\;
\sum_{\beta\in\Delta^+,\;w\to w\tau_\beta}\kern-2ex
({\omega_\alpha,\beta^\vee})\cdot X({w\tau_\beta})
\;\biggr].
\end{equation}
(Here $(\,\cdot\,,\cdot\,)$ stands for the Weyl-invariant scalar product between the characters of the maximal torus of $G$ with $\mbb R$-coefficients, $\Delta^+$ denotes the positive roots of $G$, $\omega_\alpha$ is the fundamental weight corresponding to the simple root $\alpha$, and $\beta^\vee:=2\beta/(\beta,\beta)$ is the co-root of $\beta$.)

%%%%%
%%%%%

\subsection*{Minuscule homogeneous varieties}\label{ssct:minuscule} 
Now we specialize to the case relevant for us.

\begin{m-definition}(cf. \cite[Ch. V, Proposition 2.2]{hiller-book})
A fundamental weight $\omega_\alpha$ is called \emph{minuscule} if it satisfies the condition 
\begin{equation}\label{eq:minuscule}
(\omega_\alpha,\beta^\vee)\in\{0,1\},\;\forall\,\beta\in\Delta^+.
\end{equation}
Let $P_\alpha\subset G$ be the standard parabolic subgroup defined by omitting the simple root $\alpha$. One says that $X=G/P_\alpha$ is a \emph{minuscule flag variety} if $\omega_\alpha$ is minuscule. 
\end{m-definition}
Since $P=P_\alpha\subset G$ is maximal, there is a unique Schubert divisor $\DD=\DD_\alpha$ and a unique Schubert line $\LL\cong\mbb{CP}^1$. The line bundle $\eL:=\eO_X(\DD)\in\Pic(X)$ is ample and $\Pic(X)=\mbb Z\eL$. 

The essential property of minuscule weights is that all the coefficients in Pieri's formula \eqref{eq:pieri} are either zero or one, thus 
\begin{equation}\label{eq:pieri-min}
\big[\DD\cdot X(w)\big]=\biggl[\;
\sum_{w'\in W^{P},\;w\to w'}\kern-2exX({w'})
\;\biggr].
\end{equation}
In other words, the intersection product between $\DD$ and $X(w)$ is represented by the union of the (reduced) Weyl divisors $X({w'})$, with $w\prec w'$. This is the main reason for restricting ourselves in this article to minuscule homogeneous varieties. 

\begin{m-definition}\label{def:Sd}
For $d\ges 0$, let $W^P_d:=\{w\in W^P\mid\lth(w)=\dim X-d\,\}$ and 
\begin{equation}
X_d=\underset{w\in W^P_d}{\mbox{$\bigcup$}}X(w).
\end{equation} 
be the union of $d$-dimensional Schubert subvarieties of $X$. It is a connected subvariety of $X$ (all the Schubert varieties contain the point $X(w_0^P)=X(1)$). 
For $w\in W^P_{d+1}$, consider 
\begin{equation}\label{eq:Dw}
D(w):=\underset{w\to\,w'\in W^P_d}{\mbox{$\bigcup$}}\kern-1ex X(w').
\end{equation}
Usually $D(w)$ is called \emph{the boundary} of the Schubert variety $X(w)$, and is denoted by $\partial X(w)$. Obviously, $D(w)\subset X_{d}$; the assumption that $X$ is minuscule implies  
\begin{equation}\label{eq:LXw}
\eL_{X(w)}\cong\eO_{X(w)}\big(D(w)\big). 
\end{equation}
\end{m-definition}
We will use the following important properties of Schubert varietieties:
\begin{equation}\label{eq:dw}\hspace{-20pt}
\begin{minipage}{.9\linewidth}
\begin{itemize}
\item 
$X(w),D(w)$ are (equidimensional) Cohen-Macaulay varieties \\ (cf. \cite{ram,brion-kumar}, \cite[Corollary 2.2.7]{br});
\item  
$D(w)$ is Frobenius split, actually its reduction in positive characteristics is so (cf. \textit{ibid.}).
\item 
For any $w_1,w_2\in W^P$, the scheme theoretic intersection $X(w_1)\cap X(w_2)$ is reduced (cf. \cite[Theorem 3]{ram}, \cite[Proposition 1.2.1]{brion-kumar}).
\end{itemize}
\end{minipage}
\end{equation}

\begin{m-lemma}\label{lm:iso}
{\rm(i)} 
For any $d\ges 1$, $w\in W^P_{d+1}$, the restriction homomorphisms below are bijective: 
$$
\begin{array}{r}
\mbb Z\eL=\Pic(X)\to\Pic(X(w)),\\  \mbb Z\eL=\Pic(X)\to\Pic(D(w)).
\end{array}
$$

\nit{\rm(ii)} 
For $d\ges 2$ and $w$ as above holds: 
$$
H^t(D(w),\L)=0,\;\forall\,0<t<\dim D(w)=d,\;\forall\L\in\Pic\big(D(w)\big).
$$ 
\end{m-lemma}

\begin{proof}
(i) For the first claim, the restriction homomorphism is always surjective (cf. \cite[Ch.~12]{mathieu}, \cite[Proposition 2.2.8]{br}). The injectivity is trivial: if the restriction of a line bundle to $X(w)$ is trivial, then its degree along $\LL$ vanishes. For the second claim, $D(w)$ is a union of Schubert varieties: pairwise, their intersections are reduced Schubert varieties of dimension $d-1\ges 1$ (or empty), all containing $\LL$. 
\smallskip

\nit(ii) 
We have seen that $\Pic(D(w))=\mbb Z\eL_{D(w)}$. For positive multiples of $\eL_{D(w)}$, the cohomology vanishing is \cite[Theorem 2(a)]{ram}. For large negative multiples, use that $D(w)$ is an equidimensional Cohen-Macaulay variety, so Serre duality holds. The vanishing for arbitrary (strictly) negative multiples of $\eL_{D(w)}$ is implied now by the fact that $D(w)$ is Frobenius-split (cf. \cite[Lemma 1.2.7]{brion-kumar}).
\end{proof}

%%%%%%%%%%%%%%%%%%%%%%%%%%%%%%%%%%%%%%%%%%%%
%%%%%%%%%%%%%%%%%%%%%%%%%%%%%%%%%%%%%%%%%%%%

\section{The main result}\label{sct:result}

Henceforth, for any locally free sheaf $\crl V$, we denote $\crl E:=\cEnd(\crl V)$. We are going to use the following elementary lemma:
\begin{m-lemma}\label{lm:split}
Let $S$ be a closed subscheme of a scheme $S\,'$, such that  $\Gamma(\eO_S)=\Gamma(\eO_{S'})=\bbC$. Let $\crl V_{S'}$ be a locally free sheaf on $S\,'$; we assume that the restriction $\Gamma(\crl E_{S'})\to\Gamma(\crl E_{S})$ is surjective. Then $\crl V_{S'}$ splits if and only if $\crl V_S$ splits.
\end{m-lemma} 

\begin{proof}
See \cite[\S1]{hal}.
\end{proof}

\begin{m-theorem}\label{thm:homog}
An arbitrary vector bundle $\crl V$ on $X$ splits if and only if $\crl V_{X_2}$ splits. 
\end{m-theorem}

\begin{proof}
The proof is by induction: we show that the splitting of $\crl V$ along $X_{d}$, with $d\ges 2$, implies its splitting along $X_{d+1}$. 
\smallskip

\nit\textit{Claim 1}\quad 
For $d\ges 1$, $w\in W^P_{d+1}$ holds: 
$$
H^1\big(\,X(w),\crl E\otimes\eL^{-1}\,\big)=0.
$$ 
Indeed, since $X(w)$ is Cohen-Macaulay (cf. \eqref{eq:dw}), the Serre duality implies 
$$
a_0:=\min\{a\ges 1\,|\,H^1(X(w),\crl E\otimes\eL^{-a}){=}\,0\}<+\infty.
$$
We claim that $a_0=1$. Otherwise, if $a_0\ges 2$, the exact sequence 
$$
0\to \eL_{X(w)}^{-a_0}\to \eL_{X(w)}^{-(a_0-1)}\to \eL_{D(w)}^{-(a_0-1)}\to0
$$
tensored by $\crl E$ yields: 
$\;H^1\big(\crl E\otimes\eL_{X(w)}^{-a_0}\big)\to 
H^1\big(\crl E\otimes\eL_{X(w)}^{-(a_0-1)}\big)\to 
H^1\big(\crl E\otimes\eL_{D(w)}^{-(a_0-1)}\big).$

The left-hand-side vanishes, by definition. The induction hypothesis says that $\crl E_{D(w)}$ is a direct sum of line bundles, so lemma \ref{lm:iso}(ii) implies that the right-hand-side vanishes too. Therefore the middle term vanishes, with $a_0-1\ges 1$, which contradicts the minimality of $a_0$. 
\smallskip 

\nit\textit{Claim 2}\quad 
$\crl V_{X(w)}$ splits, for all $w\in W^P_{d+1}$. 

\nit Indeed, the previous claim and the exact sequence 
$$
0\to\big(\crl E\otimes\eL^{-1}\big)_{X(w)}\to\crl E_{X(w)}\to\crl E_{D(w)}\to0
$$
yields that $\res^{X(w)}_{D(w)}:\Gamma(X(w),\crl E)\to\Gamma(D(w),\crl E)$ is surjective, and we apply the lemma \ref{lm:split}. 
\smallskip 

\nit\textit{Claim 3}\quad $\crl V_{X_{d+1}}$ splits.

\nit The splitting of $\crl V_{X_d}$ implies that there is $\vphi_d\in\Gamma(X_d,\crl E)$ with $\rk(\crl V)$ pairwise distinct eigenvalues, whose eigenspaces are the direct summands of $\crl V_{X_d}$. As $D(w)\subset X_d$, for any $w\in W^P_{d+1}$, the surjectivity of $\res^{X(w)}_{D(w)}$ implies that ${\vphi_d}|_{D(w)}$ extends to some $\vphi_d(w)\in\Gamma(X(w), \crl E)$ which still has the same eigenvalues. Now define $\vphi_{d+1}\in\Gamma(X_{d+1},\crl E)$ by setting 
$$
{\vphi_{d+1}|}_{X(w)}:=\vphi_d(w),\;\forall\,w\in W^P_{d+1}.
$$ 
We claim that $\vphi_{d+1}$ is well-defined: indeed, the intersection of any two distinct components $X(w_1), X(w_2)\subset X_{d+1}$ is reduced (cf. \eqref{eq:dw}) and 
$$ 
X(w_1)\cap X(w_2)\subset D(w_1)\cap D(w_2)\subset X_d.
$$
Moreover, we have $\vphi_d(w_j)|_{D(w_j)}=\vphi_d|_{D(w_j)}$ and $\vphi_d$ is defined over $X_d\supset D(w_1)\cap D(w_2)$.
\end{proof}

\begin{m-remark}\label{rmk:optimal}
There are stable vector bundles of arbitrary rank on $\mbb P^3$ with trivial splitting type along the generic line (the ADHM construction). Thus is not possible to verify the splitting of a vector bundle by restricting it to any finite union of test lines.
\end{m-remark}

%%%%%%%%%%%%%%%%%%%%%%%%%%%%%%%%%%%%%%%%%%%%
%%%%%%%%%%%%%%%%%%%%%%%%%%%%%%%%%%%%%%%%%%%%

\section{The list}\label{sct:expl}

In this section we perform a case-by-case analysis of the theorem \ref{thm:homog}. The minuscule varieties with cyclic Picard group are listed e.g. in \cite[\S5.2]{hiller-book} and \cite[\S2.3]{lms}. Some varieties are minuscule for several groups (e.g. $\mbb P^{2n-1}$ is minuscule for $G=\Sl(2n)$ and $\Sp(n)$). Since we are primarily interested in the variety itself, we choose the easier description ($\Sl(2n)$ in the example). As we will see, it turns out that the $2$-dimensional Schubert varieties are very simple geometric objects, although their embedding in the ambient homogeneous space can be subtle: they are mostly projective $2$-planes. 

\begin{enumerate}
\item 
$X=\Grs(k;n)$, the Grassmannian of $k$-planes in $\mbb C^n$, is minuscule for $G=\Sl(n)$. Let $e_1,\ldots,e_n$ be the standard base in $\mbb C^n$. 
\begin{enumerate}
	\item 
$k=1,n-1\;\Rightarrow\;X\cong \mbb{CP}^{n-1}\;\Rightarrow\;X_2\cong\mbb{CP}^2$. 
	\item 
$1<k<n-1\;\Rightarrow\;X_2\cong\mbb{CP}^2\cup_\LL\mbb{CP}^2$. \\ 
So, $X_2$ is the union (wedge) of the Schubert $2$-planes 
$$
\begin{array}{l}
\{\lran{e_1,\ldots,e_{k-2},f_{k-1},f_{k}}\mid 
f_{k-1},f_{k}\in\lran{e_{k-1},e_{k},e_{k+1}}\},
\\ 
\{\lran{e_1,\ldots,e_{k-2},e_{k-1},f_{k}}\mid 
f_{k}\in\lran{e_{k},e_{k+1},e_{k+2}}\}
\end{array}
$$ 
glued along the line 
$\LL=\{\lran{e_1,\ldots,e_{k-2},e_{k-1},f_{k}}\mid f_{k}\in\lran{e_{k},e_{k+1}}\}$.
\end{enumerate}
\smallskip 
%%%
\item 
$X=\euf S_n$, the Grassmannian of maximal isotropic (Lagrangian) $n$-planes in $\mbb C^{2n}$ with respect to a non-degenerate symmetric bilinear form (known as the  \emph{spinor variety}), is minuscule for $G=\SO(2n)$; its dimension equals $\frac{n(n-1)}{2}$. 

Let $\mbb C\subset\ldots\mbb C^{n-1}\subset\mbb C^n$ be a flag of isotropic subspaces in $\mbb C^{2n}$. The  unique $2$-dimensional Schubert subvariety of $X$ is: 
$$
X_2=\big\{
U\subset\mbb C^{2n}\mid \dim U=n,\;\mbb C^{n-3}\subset U,\;\dim(U\cap\mbb C^{n})\ges n-2
\big\}\cong\mbb P^2.
$$
\smallskip 
%%%
\item 
The $2n$-dimensional quadric 
$
\euf Q_{2n}:=\{x_0y_0+x_1y_1+\ldots x_ny_n=0\}\subset\mbb{CP}^{2n+1},
$ 
is minuscule for $G=\SO(2n+2)$. For $n\ges3$, the $2$-dimensional Schubert variety is 
$$
X_2=\{x_0=\ldots=x_n=y_3=\ldots=y_{n}=0\}\cong\mbb{CP}^2.
$$
For $n=2$, we have 
$X_2=\{x_0=x_1=x_2=0\}\cup\{x_0=x_1=y_2=0\}\cong\mbb{CP}^2\cup_{\mbb{CP}^1}\mbb{CP}^2$. 
\smallskip 
%%%
\item 
The Cayley plane $X=\mbb{OP}^2$ is studied in detail in \cite{iliev+manivel,robl}. It is a $16$-dimensional variety, which can be identified with the projective plane over the octonions, and is minuscule for the exceptional group $G=E_6$. We are interested only in the fact that there is only one $2$-dimensional Schubert cell, which is isomorphic to $\mbb{CP}^2$; the embedding is explicitly described in \cite[pp.~151]{iliev+manivel}.  
\smallskip 
%%%
\item 
The Freundenthal variety is a $27$-dimensional variety, minuscule for $G=E_7$. The Schubert varieties and the corresponding Hasse diagram can be found in \cite[\S3.1]{robl}. Again, there is only one \emph{smooth}, $2$-dimensional Schubert variety $X_2$. Since $X_2$ is rational and its Picard group is cyclic, we deduce that $X_2\cong\mbb{CP}^2$. 
\end{enumerate}

The conclusions of our discussion are summarized in the table below: theorem \ref{thm:homog} yields the following splitting criteria. \vskip2ex  

\begin{longtable}%
{|>{\centering}m{14ex}|%
>{\centering}m{13.5ex}|%
>{\centering}m{0.1\linewidth}|%
>{\raggedright}m{0.4\linewidth}|ccc}
\cline{1-4} 
Variety $X=G/P$ & 
group $G$ & 
$X_2\cong\ldots$ & 
\qquad splitting criterion for $\crl V$ on $X$ 
&&&
\tabularnewline \cline{1-4} 
$\mbb P^{n-1}$ & 
$\Sl(n),\,n\ges 3$ & 
$\mbb P^2$ & 
$\crl V\,\Leftrightarrow\,\crl V_{\mbb P^2}$ splits (Horrocks' criterion) 
&&&
\tabularnewline \cline{1-4} 
$\Grs(k;n)$, $1<k<n-1$ & 
$\Sl(n),\, n\ges 4$ & 
$\mbb P^2\cup_{\mbb P^1}\mbb P^2$ &
$\crl V\,\Leftrightarrow\,\crl V_{\mbb P^2\cup_{\mbb P^1}\mbb P^2}$ splits
&&&
\tabularnewline \cline{1-4} 
spinor variety $\euf S_n$&
$\SO(2n)$, $n\ges 3$&
$\mbb P^2$&
$\crl V\,\Leftrightarrow\,\crl V_{\mbb P^2}$ splits
&&&
\tabularnewline \cline{1-4}  
quadric $\euf Q_{4}$ &
$\SO(6)$&
$\mbb P^2\cup_{\mbb P^1}\mbb P^2$&
$\crl V\,\Leftrightarrow\,\crl V_{\mbb P^2\cup_{\mbb P^1}\mbb P^2}$ splits
&&
\begin{minipage}[c]{4ex}\vspace{-2.5ex}
\begin{equation}\label{eq:expl}\;\end{equation}
\end{minipage}
&
\tabularnewline \cline{1-4}  
quadric \\ $\euf Q_{2n}$, $n\ges 3$ &
$\SO(2n+2)$&
$\mbb P^2$&
$\crl V\,\Leftrightarrow\,\crl V_{\mbb P^2}$ splits
&&&
\tabularnewline \cline{1-4}  
Cayley plane $\mbb{OP}^2$&
$E_6$&
$\mbb P^2$&
$\crl V\,\Leftrightarrow\,\crl V_{\mbb P^2}$ splits
&&&
\tabularnewline \cline{1-4} 
Freundenthal variety&
$E_7$&
$\mbb P^2$&
$\crl V\,\Leftrightarrow\,\crl V_{\mbb P^2}$ splits
&&&
\tabularnewline \cline{1-4}
\end{longtable}
\vskip2ex

Let us observe that this analysis yields splitting criteria for vector bundles on the odd-dimensional quadric $\euf Q_{2n+1}$, $n\ges 2$, as well: $\euf Q_{2n}\subset\euf Q_{2n+1}$ is the hyperplane section, and $\crl V$ on $\euf Q_{2n+1}$ splits if and only if $\crl V_{\euf Q_{2n}}$ splits (cf. \cite[\S2]{hal}). We deduce that: 
\begin{equation}\label{eq:Q-odd}
\begin{array}{rcl}
\crl V\text{ on }\euf Q_{5}\text{ splits}
&\;\Leftrightarrow\;&
\crl V_{\mbb P^2\cup_{\mbb P^1}\mbb P^2}\text{ splits},
\\ 
\crl V\text{ on }\euf Q_{2n+1}, n\ges3,\text{ splits}
&\;\Leftrightarrow\;&
\crl V_{\mbb P^2}\text{ splits}.
\end{array}
\end{equation}

\subsection*{How to decide the splitting along $X_2$?} 
As the table above shows, the two dimensional Schubert subvarieties are mostly projective planes. So, one would like to have a tool to check the splitting of a vector bundle on $\mbb{CP}^2$. 

\begin{m-proposition}\label{prop:x2}
A vector bundle $\crl V$ on $\mbb P^2$ splits if and only if $H^1(\crl E(-1))=0.$ 
\end{m-proposition}

\begin{proof}
The condition is necessary, so let us prove that it is also sufficient. Let $\LL\subset\mbb P^2$ be a straight line. Then $\crl V_\LL$ decomposes into a direct sum of line bundles, by Grothendieck's theorem: thus there is $\vphi_\LL\in\Gamma(\crl E_\LL)$ with $\rk(\crl V)$ distinct eigenvalues. The hypothesis implies that $\Gamma(\crl E)\to\Gamma(\crl E_\LL)$ is surjective, so $\vphi_\LL$ extends to $\vphi\in\Gamma(\crl E)$ with the same eigenvalues.
\end{proof}
\nit Note that, in the situation \eqref{eq:expl}, the natural choice for $\LL$ is the (unique) Schubert line in $X$. For the cases when $X_2\cong\mbb{CP}^2\cup_\LL\mbb{CP}^2$, we apply the proposition for both components.

%%%%%%%%%%%%%%%%%%%%%%%%%%%%%%%%%%%%%%%%%%%%
%%%%%%%%%%%%%%%%%%%%%%%%%%%%%%%%%%%%%%%%%%%%

\section{Concluding remarks}\label{sect:end}

It is natural to ask if there is a similar splitting criterion for vector bundles on homogeneous varieties which are not necessarily minuscule. This hypothesis was used (cf. \eqref{eq:pieri-min}) to represent the class of the intersection product $\DD\cdot X(w)$ as a sum of \emph{reduced} boundary divisors of $X(w)$. 
In general, Pieri's formula \eqref{eq:pieri} involves multiplicities. For this reason, one should consider in theorem \ref{thm:homog} a union of thickenings of the Schubert surfaces, instead of the set-theoretical union $X_2$. However, the resulting form does not seem appealing. 

Let us remark that in \eqref{eq:expl} and \eqref{eq:Q-odd} there is one notable absence from the list of \emph{co-minuscule} (also called Hermitian symmetric) varieties: the Lagrangian Grassmannian $\LG(n)$. If $\omega$ is a non-degenerate, skew-symmetric form on $\mbb{C}^{2n}$, then define 
\begin{equation}\label{eq:LGn}
\LG(n):=\big\{U\subset\mbb{C}^{2n}\mid\dim(U)=n,
\;{\omega|}_U=0\;\text{(that is, $U$ is isotropic for $\omega$)}\big\}.
\end{equation}
It is homogeneous for the action of the symplectic group $\Sp(n)$, of dimension $\frac{n(n+1)}{2}$, and one may wonder whether a similar splitting criterion holds for $\LG(n)$. The same approach as in section \ref{sct:result} fails; as explained above, one has to take into account multiplicities (for $\LG(n)$, some of the coefficients in \eqref{eq:pieri} equal $2$). However, one can can still prove the following weaker result.

\begin{m-theorem}\label{thm:LGn}
A vector bundle $\crl V$ on $\LG(n)$ splits if and only if its restriction to a \emph{very generally} embedded $\LG(2)\subset\LG(n)$ does so. 
\end{m-theorem}
Note that $\euf Q_3\cong\LG(2)\subset\Grs(2;4)\cong\euf Q_4$. 
Although the statement of the theorem is in the same vein as \ref{thm:homog}, the method of its proof is totally different and considerably more difficult (cf. \cite[\S7.2]{hal2}). 

%%%%%%%%%%%%%%%%%%%%%%%%%%%%%%%%%%%%%%%%%%%%
%%%%%%%%%%%%%%%%%%%%%%%%%%%%%%%%%%%%%%%%%%%%


\begin{thebibliography}{ooo}

\bibitem{bal} \textrm{E.~Ballico}, \textit{Uniform vector bundles on quadrics}. Ann. Univ. Ferrara Sez. VII \textbf{27} (1982),135--146.

\bibitem{bgg} \textrm{I.N.~Bernstein, I.M.~Gelfand, S.I.~Gelfand}, 
\textit{Schubert cells and cohomology of the spaces G/P}. 
Russ. Math. Surv. \textbf{28} (1973), 1--26.

\bibitem{br} \textrm{M.~Brion}, \textit{Lectures on the geometry of flag varieties}. In `Topics in Cohomological Studies of Algebraic Varieties' (P.~Pragacz (ed.)), Birkh\"auser Verlag, 33--86 (2005). 

\bibitem{brion-kumar} M.~Brion, S.~Kumar, \textit{Frobenius splitting methods in geometry and representation theory}. Progress Math. \textbf{231}, Birkh\"auser (2005).

\bibitem{dem} \textrm{M.~Demazure}, 
\textit{Invariants symétriques entiers des groupes de Weyl et torsion}. 
Invent. Math. \textbf{21} (1973), 287--301.

\bibitem{ehs} \textrm{G.~Elencwajg, A.~Hirschowitz, M.~Schneider}, \textit{Les fibr\'es uniformes de rang au plus $n$ sur $\mbb P^n(\mbb C)$ sont ceux qu'on croit}. In `Vector bundles and differential equations' (Proc. Conf. Nice, 1979), pp. 37--63, Birkh\"auser, (1980).

\bibitem{guy} \textrm{M.~Guyot}, \textit{Caract\'erisation par l'uniformit\'e des fibr\'es universels sur la Grassmannienne}. Math. Ann. \textbf{270} (1985), 47--62.

\bibitem{hal2} \textrm{M.~Halic}, 
\textit{Vector bundles on projective varieties which split along q-ample subvarieties}. preprint http://arxiv.org/abs/1410.1207.

\bibitem{hal} \textrm{M.~Halic}, 
\textit{Splitting criteria for vector bundles induced by restrictions to divisors}. 
preprint http://arxiv.org/abs/1501.07101.

\bibitem{horr} \textrm{G.~Horrocks},
\textit{Vector bundles on the punctured spectrum of a local ring}. 
Proc. Lond. Math. Soc., \textbf{14} (1964), 689--713.

\bibitem{hiller-book} \textrm{H.~Hiller}, 
\textit{Geometry of Coxeter groups}. 
Research Notes in Mathematics \textbf{54}, Pitman Advanced Publishing Program (1982).

\bibitem{iliev+manivel} \textrm{A.~Iliev, L.~Manivel}, 
\textit{The Chow ring of the Cayley plane}. 
Compositio Math. \textbf{141} (2005), 146--160.

\bibitem{lms} \textrm{V.~Lakshmibai, C.~Musili, C.S.~Seshadri}, 
\textit{Geometry of $G/P\,$--~{\rm III}: Standard monomial theory for a quasi-minuscule $P$}. 
Proc. Indian Acad. Sci., Sect. A, \textbf{88} (1979), 93--177.

\bibitem{mals} \textrm{F.~Malaspina}, 
\textit{A few splitting criteria for vector bundles}. 
Ric. Mat. 57 (2008), no. 1, 55--64.

\bibitem{mathieu} \textrm{O.~Mathieu}, 
\textit{Formules de caract\`eres pour les alg\`ebres de Kac-Moody g\'en\'erales}. 
Ast\'erisque \textbf{159-160} (1988). 

\bibitem{mos} \textrm{R.~Mu\~noz, G.~Occhetta, L.~Sol\'a Conde}, 
\textit{Uniform vector bundles on Fano manifolds and applications}. 
J.~Reine Angew. Math. \textbf{664} (2012), 141--162.

\bibitem{ottv} \textrm{G.~Ottaviani}, 
\textit{Some extensions of Horrocks criterion to vector bundles 
on Grassmannians and quadrics}. Ann. Mat. Pura Appl. 155 (1989), 317--341. 

\bibitem{ram} A.~Ramanathan, 
\textit{Schubert varieties are Cohen-Macaulay}. 
Invent. Math. \textbf{80} (1985), 283--294.

\bibitem{robl} \textrm{C.~Robles}, 
\textit{Singular loci of cominuscule Schubert varieties}. 
J.~Pure Appl. Algebra \textbf{218} (2014), 745--759.

\bibitem{vdv} \textrm{A.~Van de Ven}, \textit{On uniform vector bundles}. Math. Ann. \textbf{195} (1972), 245--248.
\end{thebibliography}
\end{document}